\documentclass{article}

\usepackage[english]{babel}
\usepackage{amsmath, amssymb}
\usepackage{amsthm}
\usepackage{amsfonts}
\usepackage{enumerate}
\usepackage[latin1]{inputenc}
\usepackage{graphics}

\newtheorem{thm}{Theorem}[section]
\newtheorem{cor}[thm]{Corollary}

\newtheorem{lemm}[thm]{Lemma}
\newtheorem{defin}[thm]{Definition}
\theoremstyle{definition}
\newtheorem{rem}[thm]{Remark}

\title{Shape sensitivity analysis of the eigenvalues of the Reissner-Mindlin system\footnote{To appear in {\it SIAM Journal on Mathematical Analysis}} }
\author{Davide Buoso and Pier Domenico Lamberti}
\date{}

\begin{document}
\maketitle
\pagestyle{plain}

\noindent
{\bf Abstract:} We consider the eigenvalue problem for the Reissner-Mindlin system arising in the study of the free vibration modes of an elastic clamped plate. We provide 
quantitative estimates for the variation of the eigenvalues upon variation of the shape of the plate. We also prove analyticity results  and establish Hadamard-type formulas.
Finally, we address the problem of minimization of the eigenvalues in the case of isovolumetric domain perturbations. In the spirit of the Rayleigh conjecture   for the biharmonic operator,  we prove that balls are critical points with volume constraint for
all simple eigenvalues and the elementary symmetric functions of multiple eigenvalues.

\vspace{11pt}

\vspace{6pt}


\section{Introduction}
Let  $\Omega $ be a bounded open set in ${\mathbb{R}}^N$ with $N\ge 2$, and $t, \lambda, \mu , k>0$ be fixed parameters.  We consider the following eigenvalue
problem 
	\begin{equation}
	\label{class}
	\left\{
	\begin{array}{l l}
	-\frac{\mu}{12}\Delta\beta-\frac{\mu+\lambda}{12}\nabla\mathrm{div}\beta  -\frac{\mu k}{t^2}(\nabla w-\beta)
		= \frac{\gamma t^2}{12}\beta, & \mathrm{in}\ \Omega , \vspace{0,2cm}\\
	-\frac{\mu k}{t^2}(\Delta w -\mathrm{div}\beta)=\gamma w, & \mathrm{in}\ \Omega,\vspace{0,2cm}\\
	\beta = 0, \ \ w=0, & \mathrm{on}\ \partial\Omega ,
	\end{array}
	\right.
	\end{equation}
in the unknowns $(\beta , w)=(\beta_1, \dots , \beta_N , w)$ (the eigenvector) and $\gamma $ (the eigenvalue).    According to the Reissner-Mindlin model for moderately thin plates, for $N=2$ system (\ref{class}) describes the free vibration modes of an elastic clamped
plate $\Omega \times (-t/2,t/2)$ with midplane $\Omega$ and thickness $t$. In that  case  
  $\lambda$ and $\mu$ are the Lam\'{e} constants, $k$ is the correction factor, $w$ the transverse displacement of the midplane, $\beta =(\beta_1,\beta_2)$   
the fiber rotation and $\gamma t^2$ the vibration frequency. We refer to Dur\'{a}n et al.~\cite{durnie} for more information and references, see also Hervella-Nieto~\cite{herthesis}.  Although $N=2$ seems to be the case of main interest in applications, our methods allow us to treat the general case without  any restriction on the space dimension.

It is well-known that the spectrum of the Reissner-Mindlin system is discrete,   hence problem (\ref{class}) has a divergent sequence of positive eigenvalues of finite multiplicity
$$
0< \gamma_{1,t}[\Omega ]\le  \gamma_{2,t}[\Omega ]\le \dots \le \gamma_{n,t}[\Omega ]\le \dots 
$$
depending on $t$ and $\Omega$.  Here each eigenvalue is repeated according to its multiplicity. 

The behavior of the solutions to Reissner-Mindlin systems as $t\to 0$ is well-known. We refer to  Brezzi and Fortin~\cite{brezzi1, brezzi2}  for a deep analysis of related 
computational problems and references. See also Lovadina et al.~\cite{lova}. In particular, it is proved  in  \cite{durnie} for $N=2$ that   $\gamma_{n,t}[\Omega]\to \gamma_{n,0}[\Omega]$ as $t\to 0$, where
$\gamma_{n,0}[\Omega]$ are the eigenvalues of the  problem 
\begin{equation}
	\label{classbiha}
	\left\{
	\begin{array}{l l}
	\frac{2\mu +\lambda}{12}\Delta^2w=\gamma w, & \mathrm{in}\ \Omega,\vspace{0,2cm}\\
	w=\nabla w =0 & \mathrm{on}\ \partial\Omega .
	\end{array}
	\right.
\end{equation}

In this paper we are interested in the dependence of $\gamma_{n,t}[\Omega]$ on $\Omega$.  
In Section~3, we provide stability estimates in the spirit of \cite{buda, bula, buladir, bulahigh, bulahighsh}. These estimates  allow to control the variation of $\gamma_{n,t}[\Omega]$ upon variation of $\Omega$.

First, we consider the case of domain deformations of the form $\phi (\Omega )$ where $\phi $ is a diffeomorphism  of class $C^{1,1}$ and in Theorem~\ref{loclip} we prove the existence of a constant $c>0$ independent of $n$ and $t$ 
such that 
\begin{equation}\label{generalest}
|\gamma_{n,t}[\phi (\Omega)]-\gamma_{n,t}[\Omega]|\le c \gamma_{n,t}[\Omega] \delta (\phi ),
\end{equation}
provided $\delta (\phi )< c^{-1}$, where $\delta (\phi )$ is defined by
\begin{equation}\label{deltafi}
\delta (\phi )= \max_{1\le |\alpha |\le 2}\sup_{x\in \Omega }|D^{\alpha}( \phi (x)- x)|.
\end{equation}

Second, we prove estimates in terms of explicit geometric quantities which measure the vicinity of two open sets $\Omega_1$ and $\Omega_2$. To do so, we assume that 
$\Omega_1$ and $\Omega_2$ belong to the same uniform class $C({\mathcal{A}})$ where ${\mathcal{A}}$ is a fixed atlas  by the help of which the open sets 
are described locally as the subgraphs of suitable continuous functions, see Definition~\ref{atclass}. In this case, it is possible to prove the existence of a constant $c>0$ independent of  $n$ and $t$ such that 
 \begin{equation}
\label{dirthm1}
| \gamma_{n,t}[\Omega_1]-\gamma_{n  ,t}[\Omega_2] |\le c \max\{ \gamma_{n,t}[\Omega_1] , \gamma_{n  ,t}[\Omega_2]\} d_{\mathcal{A}}(\Omega_1,\Omega_2),
\end{equation}
provided $d_{\mathcal{A}}(\Omega_1,\Omega_2)\le c^{-1}$, where $d_{\mathcal{A}}(\Omega_1,\Omega_2)$ is the so-called atlas distance of $\Omega_1$ and $\Omega_2$. See Theorem~\ref{dirthm} and Definition~\ref{dev}.
We note that the  atlas distance $d_{\mathcal{A}}(\Omega_1,\Omega_2)$ is an easily computable one-dimensional distance which measures the gap between the graphs describing the boundaries of $\Omega_1,\Omega_2$ and that it is possible to control it via the more familiar Hausdorff distance between $\partial \Omega_1$ and $\partial \Omega_2$, see Theorem~\ref{metricsurv}.     Importantly,  the atlas class $C({\mathcal{A}})$ includes open sets with strong boundary degenerations such as cusps of exponential type. In fact, if the modulus of continuity $\omega $ of the
functions describing the boundaries of $\Omega_1$ and $\Omega_2$ is fixed and the boundary of one of the two domains is contained in an $\epsilon$-neighborhood of the boundary of the other one,   then it is possible to prove an estimate via $\omega (\epsilon )$. See Corollary~\ref{inclestthm}.

Note that  error estimates independent of $t$ for a finite element discretization of the eigenvalue problem (\ref{class}) on a polygon in the plane have been obtained in Dur\'{a}n et al.~\cite{durnie}. Considering that polygons are typically used in order to approximate sufficiently regular  planar domains, we believe that our estimates  complement those in 
\cite{durnie}.

In Section 4,  we consider families of open sets $\phi (\Omega)$ parametrized by Lipschitz homeomorphisms $\phi$, and we prove analyticity results  for the dependence of $\gamma_{n,t}[\phi (\Omega )]$ on $\phi$. Following the analysis of \cite{lala}, we prove that simple eigenvalues and the elementary symmetric functions of multiple eigenvalues depend real analytically on $\phi$, and we establish Hadamard-type formulas for the Fr\'{e}chet differentials.  See Theorem~\ref{thesame}. In particular, if $\Omega$ is sufficiently smooth and  $\gamma_{n,t}[ \Omega ]$ is simple then for perturbations of the identity $I$ of the type $\phi_{\epsilon}=I+\epsilon \psi$, $\epsilon\in {\mathbb{R}} $, we have
\begin{equation}\label{hada}{\frac{d  \gamma_{n,t}[\phi_{\epsilon}(\Omega) ]}{d\epsilon}}_{|_{\epsilon =0}}=
-\int_{\partial\Omega}\biggl(
		\frac{\mu}{12}\biggl|\frac{\partial\beta}{\partial n}\biggr|^2
	+\frac{\mu+\lambda}{12}\left(\frac{\partial\beta}{\partial n}\cdot  n\right)^2
	+\frac{\mu k}{t^2}\left(\frac{\partial w}{\partial n}\right)^2\biggr)\psi\cdot n d\sigma,
\end{equation}
where $n$ is the unit outer normal to $\partial\Omega$ and $(\beta , w)$ is an eigenvector associated with $\gamma_{n,t}[ \Omega ]$ normalized by the condition $\int_{\Omega }w^2+\frac{t^2}{12}|\beta|^2dx=1$. The bifurcation phenomenon which occurs in the case of multiple eigenvalues is more involved and is described by the Rellich-Nagy-type Theorem~\ref{nagy}. 

Finally, in Section 5 we address the problem of the optimization of the eigenvalues in case of isovolumetric perturbations. Recall that the celebrated Rayleigh conjecture  
states that, among all bounded domains with fixed measure, the first eigenvalue of problem (\ref{classbiha}) is minimized by the ball. Such conjecture has been proved for $N=2$ by N.S.~Nadirashvili and for $N=2,3$ by M.S.~Ashbaugh and R.D.~Benguria. We refer to Henrot~\cite{henrot} for a survey on this topic. Taking into account the limiting behavior of the Reissner-Mindlin eigenvalues as $t\to 0$, it would be natural to state the same conjecture also for the Reissner-Mindlin system.  Here we give support to it  by proving that the Reissner-Mindlin system exhibits the same symmetry property of biharmonic and polyharmonic operators, see \cite{buoso, buosoplates, lalcri}.  Namely, we prove that balls are critical points with volume constraint for all  simple eigenvalues and all  symmetric functions of multiple eigenvalues of system (\ref{class}). See Theorem~\ref{puntini}. To do so, we characterize critical open sets  as  those open sets for which a suitable overdetermined system has nontrivial solutions and we prove that such overdetermined conditions are satisfied when the open set is a ball. 


\section{Preliminaries and Notation}

In this section we introduce the eigenvalue problem under consideration and the classes of open sets   which allow us to prove the quantitative estimates of Section 3. 

\subsection{The Reissner-Mindlin eigenvalue problem}

Let $\Omega $ be an open set in ${\mathbb{R}}^N$ with finite measure. By $H^1_0(\Omega )$ we denote the closure  in the
standard Sobolev space $H^1(\Omega)$ of the space of $C^{\infty }$-functions with compact support in $\Omega$. 
We set ${\mathcal{V}}(\Omega)= (H^1_0(\Omega ))^{N}\times H^1_0(\Omega )$ and we denote by $(\beta , w)$ the generic element of ${\mathcal{V}}(\Omega)$, where  $\beta =(\beta_1,\dots , \beta_N)\in (H^1_0(\Omega ))^N$ and $w\in H^1_0(\Omega )$.   

For any fixed $t, \lambda, \mu , k>0$,  we consider the weak formulation of problem (\ref{class}). Namely, we  say that $\gamma\in {\mathbb{R}} $ is an eigenvalue of the Reissner-Mindlin system 
if and only if there exists $(\beta , w)\in {\mathcal{V}}(\Omega )$ with $(\beta , w)\ne 0$ such that 
	\begin{multline}
	\label{reis}
	\frac{\mu}{12}\int_{\Omega}\nabla\beta :\nabla\eta dx+\frac{\mu+\lambda}{12}\int_{\Omega}\mathrm{div}\beta \mathrm{div}\eta dx
		+\frac{\mu k}{t^2} \int_{\Omega}(\nabla w -\beta )\cdot (\nabla v-\eta)dx
		\\=\gamma \int_{\Omega}\left(wv+\frac{t^2}{12}\beta\cdot\eta\right)dx,
	\end{multline}
for all test functions $ (\eta , v)\in {\mathcal{V}}(\Omega)$, in which case $(\beta , w)$ is called an eigenvector associated with $\gamma$.  Here by $A: B$ we denote the Frobenius product of two matrices $A,B$, defined by $A:B=\sum_{i,j=1}^Na_{ij}b_{ij}$.  Note that $\beta$ is thought as a row vector.

As customary in Spectral Theory we interpret problem (\ref{reis}) as an eigenvalue problem for a non-negative selfadjoint operator in Hilbert space as follows.    
For any fixed $t>0$, we denote by ${\mathcal{L}}^2_t(\Omega)$ the space $ (L^2(\Omega ))^N\times L^2(\Omega)$ endowed  with the scalar product $<(\beta , w),( \eta , v)>_{\Omega ,t}$ defined by the right-hand side of (\ref{reis}) (without $\gamma$) for any 
$(\beta , w),( \eta , v)\in {\mathcal{L}}^2_t(\Omega)$.  Clearly, for each $t>0$ the norm induced by such scalar  product is equivalent to the standard $L^2$-norm. 
Moreover, we consider the bilinear form $Q_{\Omega ,t }$ defined on  ${\mathcal{V}}(\Omega)\times {\mathcal{V}}(\Omega )$ by  the left-hand side of equality (\ref{reis}). 
We also denote by  $Q_{\Omega, t}(\beta , w)= Q_{\Omega, t}((\beta , w),(\beta , w))$ the quadratic form associated with the bilinear form $Q_{\Omega, t}$ and we observe that such quadratic form is coercive in ${\mathcal{V}}(\Omega)$. In particular, the corresponding norm $Q_{\Omega ,t}^{1/2}(\cdot) $ is equivalent to the standard Sobolev norm in $ {\mathcal{V}}(\Omega)$. 
This implies that the quadratic form  $Q_{\Omega, t }(\cdot )$ is closed in ${\mathcal{L}}^2_t(\Omega )$ hence  (see e.g., Davies~\cite[Ch.~4]{daviesbook})  there exists 
a non-negative self-adjoint operator $R_{\Omega ,t}$ densely defined on  ${\mathcal{L}}^2_t(\Omega )$  such the domain ${\rm Dom }(R_{\Omega ,t}^{1/2})$  of the square root $R^{1/2}_{\Omega ,t}$ of  $R_{\Omega ,t}$ is ${\mathcal{V}}(\Omega)$
and such that $Q_{\Omega, t}((\beta ,w),(\eta , v))=<R^{1/2}_{\Omega ,t}( \beta ,w), R^{1/2}_{\Omega ,t}( \eta ,v)>_{\Omega, t} $ for all  $(\beta, w),( \eta , v)\in {\mathcal{V}}(\Omega )$. In particular, 
 $(\beta , w)\in {\rm Dom }(R_{\Omega ,t})$ if and only if $(\beta , w)\in {\rm Dom }(R_{\Omega ,t}^{1/2})$ and there exists $(\theta , f)\in  L^2(\Omega)\times (L^2(\Omega))^N$ such that 
$Q_{\Omega, t}((\beta , w),(\eta , v))=<(\theta , f), (\eta , v)  >_{\Omega , t} $  for all $(\eta , v)\in {\mathcal{V}}(\Omega )$, in which case $R_{\Omega ,t}( \beta , w)=(\theta , f)$.

It follows that the eigenvalues and the eigenvectors of problem (\ref{reis}) coincide with the eigenvalues and the eigenvectors of the operator $R_{\Omega , t}$. Moreover, since $|\Omega |$ is finite, ${\mathcal{V}}(\Omega )$ is compactly embedded into ${\mathcal{L}}^2_t(\Omega )$, hence  the spectrum of $R_{\Omega ,t}$ is discrete and consists of a divergent sequence of positive eigenvalues of
finite multiplicity, which we denote by $\gamma_{n,t}[\Omega ]$, $n\in {\mathbb{N}}$. 
We note that by the Courant Min-Max Principle, we have
\begin{equation}
\gamma_{n,t}[\Omega]=\min_{\substack{E\subset {\mathcal{V}}(\Omega)\\ {\rm dim }E=n }} \max_{(\beta , w)\in E\setminus \{0\}} \frac{ Q_{\Omega, t}(\beta , w) }{{}\ \ \| ( \beta , w )\|^2_{{\mathcal{L}}^2_t(\Omega )} }
\end{equation}
for all $n\in {\mathbb{N}}$. \\

\begin{rem}   Assume that  $\Omega $ is a bounded open set in ${\mathbb{R}}^2$ representing the midplane of an elastic clamped plate $\Omega \times ]-t/2,t/2[$ of thickness $t$. 
In the literature (cf. e.g. \cite{durnie}), the weak formulation of the eigenvalue problem for the Reissner-Mindlin system describing the free vibration modes of such plate  can be found in the form
\begin{multline}\label{duran} t^3a(\beta ,\eta )
+{\mathcal{K}} t\int _{\Omega }(\nabla w-\beta )\cdot (\nabla v-\eta)dx
=\omega^2\left(t\int_{\Omega }wvdx+\frac{t^3}{12}\int_{\Omega}\beta \cdot \eta dx\right) ,
\end{multline}
where $a(\beta ,\eta )=\frac{E}{12(1-\nu^2)}\int_{\Omega}[(1-\nu)\epsilon (\beta): \epsilon (\eta)+\nu {\rm div}\beta {\rm div}\eta ]dx$ and 
${\mathcal{K}}=\frac{Ek}{2(1+\nu)}$.  Here $\omega $ is the angular vibration frequency, $\epsilon (\beta)=(\nabla \beta+\nabla^t \beta)/2$ is the linear strain tensor, $\nu $ the Poisson ratio, $E$ the Young modulus, ${\mathcal{K}}$  the shear modulus and $k$ the correction  factor (usually $k=5/6$). 
By recalling Korn's indentity 
$$
2\int_{\Omega }\epsilon (\beta ): \epsilon (\eta)dx = \int_{\Omega }\nabla \beta : \nabla \eta dx+\int_{\Omega }{\rm div }\beta {\rm div}\eta dx,
$$
which holds for any $\beta ,\eta \in {\mathcal{V}}(\Omega)$, problem (\ref{duran}) can be easily rewritten in the form (\ref{reis}) by setting $\gamma =\omega^2/t^2$ and choosing
\begin{equation}\lambda =\frac{\nu E}{1-\nu^2},\ \ {\rm and }\ \ 
\mu = \frac{E}{2(1+\nu)}.  
\end{equation} 
The formulation in (\ref{reis}) is somewhat  more general since it allows any choice of the constants $\lambda, \mu >0$ including the standard Lam\'{e} constants 
$\lambda = \nu E/ [(1+\nu)(1-2\nu)]$, $\mu = \frac{E}{2(1+\nu)}$.

We refer also to  Bathe~\cite{bat} for further details.  

\end{rem}

\subsection{The atlas class and the atlas distance}

For any set $V$ in ${\mathbb{R}}^N$ and $\delta >0$ we denote by $V_{\delta }$ the set $\{x\in V:\ d(x, \partial V )>\delta \}$. We shall also denote by $V^{\delta }$ the set
$\{x\in {\mathbb{R}}^N:\ d(x, V)<\delta \}$. Here $d(x, A)$ denotes the Euclidean distance from $x$ to a set $A$.  We recall the following definition from \cite{bulahigh}, where by cuboid we mean a set which is the isometric image of a set of the form $\Pi_{i=1}^N]a_i,b_i[$. 

\begin{defin} {\bf (Atlas Class)}
\label{atclass} 
Let $ \rho >0$, $s,s'\in\mathbb{N}$, $s'\le s$
and  $\{V_j\}_{j=1}^s$ be a family of bounded open cuboids  and
$\{r_j\}_{j=1}^{s} $ be a family of isometries in ${\mathbb{R}}^N $. 
We say that that ${\mathcal{A}}= (  \rho , s,s', \{V_j\}_{j=1}^s, \{r_j\}_{j=1}^{s} ) $ is an atlas in ${\mathbb{R}}^N$ with the parameters
$\rho , s,s', \{V_j\}_{j=1}^s, \{r_j\}_{j=1}^{s}$, briefly an atlas in ${\mathbb{R}}^N$.

We denote by $C( {\mathcal{A}}   )$ the family of all open sets $\Omega $ in ${\mathbb{R}}^N$
satisfying the following properties:

(i) $ \Omega\subset \bigcup\limits_{j=1}^s(V_j)_{\rho}$ and $(V_j)_\rho\cap\Omega\ne\emptyset;$

(ii) $V_j\cap\partial \Omega\ne\emptyset$ for $j=1,\dots s'$, $ V_j\cap \partial\Omega =\emptyset$ for $s'<j\le s$;

(iii) for $j=1,...,s$
$$
r_j(V_j)=\{\,x\in \mathbb{R}^N:~a_{ij}<x_i<b_{ij}, \,i=1,....,N\},
$$

\noindent and

$$
r_j(\Omega\cap V_j)=\{x\in\mathbb{R}^N:~a_{Nj}<x_N<g_{j}(\bar
x),~\bar x\in W_j\},$$

\noindent where $\bar x=(x_1,...,x_{N-1})$, $W_j=\{\bar
x\in\mathbb{R}^{N-1}:~a_{ij}<x_i<b_{ij},\,i=1,...,N-1\}$
and $g_j$ is a continuous function defined on $\overline {W}_j$ (it is meant that if $s'<j\le s$ then $g_j(\bar x)=b_{Nj}$ for all $\bar x\in \overline{W}_j$);
moreover for $j=1,\dots ,s'$
$$
a_{Nj}+\rho\le g_j(\bar x)\le b_{Nj}-\rho ,$$

\noindent for all $\bar x\in \overline{W}_j$.

We say that an open set $\Omega$ in ${\mathbb{R}}^N$ is an open set with a continuous boundary if $\Omega $ is of class
$C( {\mathcal{A}}   ) $ for some atlas $ {\mathcal{A}}   $.

Let $\omega :[0,\infty [\to [0,\infty [ $ be a modulus of continuity, i.e., a  continuous non-decreasing function
such that $\omega (0)=0$ and, for some $k>0$, $\omega (t)\geq kt$ for all $0\le t\le 1$.
Let $M>0$.
We denote by  $C_M^{\omega (\cdot )}( {\mathcal{A}}   )$ the family of all open sets  $\Omega$ in ${\mathbb{R}}^N$ belonging to
$C( {\mathcal{A}}   )$ and such that all the functions $g_j$ in Definition \ref{atclass} $(iii)$ satisfy the condition
\begin{equation}
\label{omegafamcond}
|g_j(\bar x) -g_j(\bar y)|\le M \omega (| \bar x -\bar y| ),
\end{equation}
for all $\bar x, \bar y \in {\overline{W}}_j$.

We also say that an open set is of class $C^{\omega (\cdot )}$ if there exists an atlas ${\mathcal{A}}$ and $M>0$ such that $\Omega \in C^{\omega (\cdot )}_M({\mathcal{A}} )$.

\end{defin}

The  family  of open sets of class $C( {\mathcal{A}} ) $ can be thought as a metric space endowed with so-called Atlas Distance. We recall  the definition introduced in \cite{bulahigh}. 

\begin{defin} {\bf (Atlas distance)}
\label{dev}
Let ${\mathcal{A}} =(\rho , s, s', \{V_j\}_{j=1}^s , \{r_j\}_{j=1}^s )$ be an atlas in ${\mathbb{R}}^N$.
For all $\Omega_1 , \Omega_2\in C({\mathcal{A}})$ we define the `atlas distance' $d_{{\mathcal{A}}}$ by

\begin{equation}\label{dev1}
d_{{\mathcal{A}}}(\Omega_1,\Omega_2)=\max_{j=1, \dots ,s}\sup_{(\bar x , x_N)\in r_j(V_j)}\left| g_{1j}(\bar x) - g_{2j}(\bar x)  \right|,
\end{equation}
where $g_{1j}$, $g_{2j}$ respectively, are the functions describing the boundaries of $\Omega_1, \Omega_2$ respectively, as in Definition \ref{atclass} $(iii)$.
\end{defin}

The atlas distance  depends on the chosen atlas but has the advantage of being easily computable. Moreover, we observe that it can be  controlled via the Hausdorf distance. Indeed,  we have the following theorem where, for the sake of completeness, we collect also other relevant properties of the atlas distance proved in \cite{bulahigh}. 

 Given two sets  $A,B$ is $ {\mathbb{R}}^N$ the 
 lower Hausdorff-Pompeiu deviation of $A$ from $B$ is defined in \cite{bulahigh}  by $
{\mathit{d}}_{{\mathcal{H}}{\mathcal{P}}} (A,B)=\min \{ \sup_{x\in A}d(x, B),\, \sup_{x\in B}d(x, A)\}$. Note that the standard Hausdorff-Pompeiu 
distance of $A$ and $B$ is ${\mathit{d}}^{{\mathcal{H}}{\mathcal{P}}} (A,B)=\max \{ \sup_{x\in A}d(x, B),\, \sup_{x\in B}d(x, A)\}$. 

\begin{thm}\label{metricsurv} Let ${\mathcal{A}} =(\rho , s, s', \{V_j\}_{j=1}^s , \{r_j\}_{j=1}^s )$ be an atlas, $\omega $ a modulus of continuity as in Definition~\ref{atclass} and $M>0$. Let   $\tilde {\mathcal{A}}=(    \rho /2 ,$ $  s,s',$ $  \{ (V_j)_{\rho /2}\}_{j=1}^s,$ $  \{r_j\}_{j=1}^{s} )   $. Then the following statements hold:
\begin{itemize}
\item[(i)] $ (C({\mathcal{A}}), d_{{\mathcal{A}}})$ is a complete metric space;
\item[(ii)] $C_M^{\omega (\cdot )}( {\mathcal{A}}   )$ is a compact subset of $C({\mathcal{A}})$;
\item[(iii)] There exists $c >0$ depending only on $N, {\mathcal{A}} ,
 \omega ,  M$  such that
\begin{equation}
\label{hausdorff}
d^{{\mathcal{H}}   {\mathcal{P}}}(\partial \Omega_1,\partial \Omega_2)\le
d_{\tilde {\mathcal{A}}}(\Omega_1, \Omega_2) \le c\, \omega (d_{{\mathcal{H}}   {\mathcal{P}}}(\partial \Omega_1,\partial \Omega_2)),
\end{equation}
for all  $\Omega_1, \Omega_2 \in C^{\omega (\cdot ) }_M({\mathcal{A}})$.
\end{itemize}
\end{thm}

\section{Quantitative estimates}

\subsection{Estimates via diffeomorphisms}

Given an open set $\Omega$ in ${\mathbb{R}}^N$ with finite measure, we consider a diffeomorphism from $\Omega$ onto another open set $\phi (\Omega)$ in ${\mathbb{R}}^N$ and we 
prove a quantitative stability  estimate for $|\gamma_{n,t}[\phi (\Omega ) ]-\gamma_{n,t}[\Omega  ]|$ in terms of  the measure of vicinity $\delta (\phi )$ defined by (\ref{deltafi}).
In order to obtain an estimate independent of $t$, we use the special 
transformation ${C_{\phi }}$ from the space ${\mathcal{V}}(\Omega )$ onto ${\mathcal{V}}(\phi (\Omega ))$ defined by
\begin{equation}\label{cifi0}
C_{\phi }(\beta , w)=  (\beta \nabla\phi^{-1}, w )\circ \phi ^{(-1)},
\end{equation}
for all $(\beta ,w)\in {\mathcal{V}}(\Omega)$.  Here and in the sequel we denote by $A^{-1}$ the inverse of a matrix $A$, as opposed to the inverse of a function $f$ which is denoted by $f^{(-1)}$; we shall also denote by $A^T$ the transpose of  $A$.

It is clear that in order to guarantee that $C_{\phi}$ is  well-defined, it suffices to assume that $\phi$ is a diffeomorphism of class $C^{1,1}$, i.e., $\phi$ and its inverse have Lipschitz continuous gradients. In fact, it is easy to prove 
the following lemma that will be used in the sequel.

\begin{lemm}\label{cifi}  Let $\Omega$  be an open set in ${\mathbb{R}}^N$ and let $\phi:\Omega \to \phi (\Omega ) $ be a diffeomorphism of class $C^{1,1}$ from $\Omega$ onto an open set $\phi (\Omega)$ in ${\mathbb{R}}^N$. Assume that 
$$
\max_{1\le |\alpha |\le 2}\sup_{x\in \Omega }|D^{\alpha} \phi (x)|<\infty,\ \ \  \inf_{x\in \Omega }|{\rm det}\nabla \phi (x)|>0. 
$$ 
Then $C_{\phi}$ is a linear homeomorphism from ${\mathcal{V}}(\Omega )$ onto ${\mathcal{V}}(\phi (\Omega ))$.
\end{lemm}

Then we can prove the following

\begin{lemm}\label{fiest} Let $\Omega$ be an open set in ${\mathbb{R}}^N$ with finite measure and let $\phi:\Omega \to \phi (\Omega ) $ be a diffeomorphism of class $C^{1,1}$ from $\Omega$ onto an open set $\phi (\Omega)$ in ${\mathbb{R}}^N$. Assume that there exist $M_1,M_2>0$ such that 
\begin{equation}\label{fiest1}
\max_{1\le |\alpha |\le 2}\sup_{x\in \Omega }|D^{\alpha} \phi (x)|<M_1,\  \ \ \inf_{x\in \Omega }|{\rm det}\nabla \phi (x)|>M_2,
\end{equation} 
for all $x\in\Omega$. Then there exists $c>0$ depending only on $N, M_1, M_2, \lambda, \mu$ and $|\Omega |$ such that 
\begin{equation}\label{fiest2}
|Q_{\phi(\Omega), t}(C_{\phi }(\beta , w))-Q_{\Omega, t}(\beta , w)|\le c Q_{\Omega , t}(\beta , w )\delta (\phi ),
\end{equation}
for all $t>0$ and $(\beta , w)\in {\mathcal{V}}(\Omega )$.
\end{lemm}

{\bf Proof. } Let $(\beta , w)\in {\mathcal{V}}(\Omega )$.  To shorten our notation, we denote by $C_{\phi }^{(1)}(\beta )$ the first entry of $C_{\phi }(\beta , w)$, i.e., 
$C_{\phi }^{(1)}(\beta )  =  (\beta\nabla\phi^{-1}  )\circ \phi ^{(-1)}$.  We begin by estimating 
$\int_{\phi (\Omega )}|\nabla C_{\phi }^{(1)}(\beta ) |^2 dy - \int_{\Omega }|\nabla \beta  |^2 dx$. 
By means of a change of variables, we get
\begin{equation}\label{fiest4}
\int_{\phi (\Omega )}|\nabla C_{\phi }^{(1)}( \beta ) |^2 dy=\int_{\Omega }|  (\nabla (\beta\nabla\phi^{-1}))\nabla \phi ^{-1}|^2|{\rm det}\nabla \phi |dx.
\end{equation}
It is easy to see that in order to estimate $\int_{\phi (\Omega )}|\nabla C_{\phi }^{(1)}(\beta ) |^2 dy - \int_{\Omega }|\nabla \beta  |^2 dx$ it suffices to estimate 
$ \int_{\Omega }(|  (\nabla (\beta\nabla\phi^{-1}))\nabla \phi ^{-1}|^2-|\nabla \beta |^2)|{\rm det}\nabla \phi |dx $. We clearly have that 
\begin{eqnarray}\label{fiest5}\lefteqn{
| \int_{\Omega }(|  (\nabla (\beta\nabla\phi^{-1}))\nabla \phi ^{-1}|^2-|\nabla \beta |^2)|{\rm det}\nabla \phi |dx | }\\
& &\qquad\le \| {\rm det }\nabla \phi \|_{L^{\infty }(\Omega )}  \|   (\nabla (\beta\nabla\phi^{-1}))\nabla \phi ^{-1} -\nabla \beta   \|_{L^2(\Omega )} \nonumber \\ 
& &\qquad\quad \cdot  (  \|   (\nabla (\beta\nabla\phi^{-1}))\nabla \phi ^{-1}   \|_{L^2(\Omega )} +  \|   \nabla \beta   \|_{L^2(\Omega )}  ). \nonumber 
\end{eqnarray}
By the triangle inequality we get 
\begin{eqnarray}\label{fiest6}\lefteqn{  \|   (\nabla (\beta\nabla\phi^{-1}))\nabla \phi ^{-1} -\nabla \beta   \|_{L^2(\Omega )}      }\\
& & \le \| \nabla \phi ^{-1}\|_{L^{\infty }(\Omega )}\| \nabla (\beta   \nabla\phi^{-1}) -\nabla \beta      \|_{L^2(\Omega )}\nonumber \\
& & \quad + \|  \nabla\phi^{-1} -I \|_{L^{\infty }(\Omega )}\| \nabla \beta \|_{L^2(\Omega )}\nonumber  
\end{eqnarray}
and 
\begin{eqnarray}\label{fiest7}\| \nabla (\beta   \nabla\phi^{-1}) -\nabla\beta     \|_{L^2(\Omega )}& \le &  \|   \nabla\phi^{-1} -I \|_{L^{\infty }(\Omega )}\| \nabla \beta \|_{L^2(\Omega )}     \\
& +&   \|\nabla (  \nabla\phi^{-1}  )\|_{L^{\infty }(\Omega )}\|  \beta \|_{L^2(\Omega )} .     \nonumber  
\end{eqnarray}
Moreover
\begin{eqnarray}\label{fiest8}
 \|   \nabla (\beta\nabla\phi^{-1})  \|_{L^2(\Omega )} &  \le &  \|   \nabla\phi^{-1}  \|_{L^{\infty }(\Omega )}\| \nabla \beta \|_{L^2(\Omega )}  \\   
 & + &   \|\nabla (   \nabla\phi^{-1})  \|_{L^{\infty }(\Omega )}\|  \beta \|_{L^2(\Omega )} . \nonumber 
\end{eqnarray}

By standard calculus it follows that there exists a constant $c>0$ depending only on $N, M_1, M_2$ such that 
\begin{equation}\label{fiest9}
\| \nabla \phi ^{-1}\|_{L^{\infty }(\Omega )}\le c
\end{equation}
and 
\begin{equation}\label{fiest10}
\|   \nabla\phi^{-1} -I \|_{L^{\infty }(\Omega )},\   \|\nabla (  \nabla\phi^{-1})  )\|_{L^{\infty }(\Omega )} \le c\delta (\phi).
\end{equation}

By using the Poincar\'{e} inequality $\|  \beta \|_{L^2(\Omega )}\le c \|  \nabla \beta \|_{L^2(\Omega )}$ with $c$ depending only on $N$  and $|\Omega|$, and combining inequalities (\ref{fiest4})-(\ref{fiest10}) we conclude that 
\begin{eqnarray}
\label{fiest11}\biggl|\int_{\phi (\Omega )}|\nabla C_{\phi }^{(1)}(\beta ) |^2 dy - \int_{\Omega }|\nabla \beta  |^2 dx \biggr|\le c_1\delta (\phi )\int_{\Omega }|\nabla \beta |^2dx,
\end{eqnarray}
where the constant $c_1$ depends only on $N, M_1, M_2$ and $|\Omega |$. 

Similarly, one can also prove the existence of a constant $c_2>0$ depending only on $N, M_1, M_2$ and $|\Omega |$ such that 
\begin{eqnarray}
\label{fiest12}\biggl|\int_{\phi (\Omega )}({\rm div}\, C_{\phi }^{(1)}(\beta ))^2  dy - \int_{\Omega }({\rm div} \beta)^2  dx \biggr|\le c_2\delta (\phi )\int_{\Omega }|\nabla \beta |^2dx.
\end{eqnarray}

Finally, we estimate 
$
\int_{\phi (\Omega )}| \nabla (w\circ \phi^{(-1)}) -   C_{\phi }^{(1)}(\beta )  |^2dy - \int_{\Omega }|\nabla w -\beta |^2dx. 
$
We note that 
\begin{equation}\label{fiest13}
\int_{\phi (\Omega )}| \nabla (w\circ \phi^{(-1)}) -   C_{\phi }^{(1)}(\beta )  |^2dy=\int_{\Omega }|(\nabla w-\beta )\cdot \nabla \phi ^{-1}  |^2 |{\rm det}\nabla \phi |dx
\end{equation}
and that 
\begin{eqnarray}\label{fiest14}
\lefteqn{\int_{\Omega }|\, |(\nabla w-\beta )\cdot \nabla \phi ^{-1}  |^2 -|\nabla w -\beta |^2|dx}\nonumber  \\
& & \qquad\qquad\qquad \le \|  \nabla \phi ^{-1}(\nabla \phi^{-1} )^{T} -I \|_{L^{\infty }(\Omega )}\int_{\Omega }|\nabla w -\beta |^2dx.
\end{eqnarray}
It follows that there exists $c_3>0$ depending only on $N, M_1, M_2$ such that 
\begin{equation}\label{fiest15}
\biggl|\int_{\phi (\Omega )}| \nabla (w\circ \phi^{(-1)}) -   C_{\phi }^{(1)}(\beta )  |^2dy - \int_{\Omega }|\nabla w -\beta |^2dx\biggr|\le c_3\delta (\phi ) \int_{\Omega }|\nabla w -\beta |^2dx.
\end{equation}
By combining inequalities (\ref{fiest11}), (\ref{fiest12}), (\ref{fiest15}), we deduce the validity of (\ref{fiest2}). \hfill $\Box$\\

As in the case of elliptic partial differential equations discussed in \cite{bulahigh}, we can prove the following 

\begin{thm}\label{loclip} Let $\Omega $ be an open set in ${\mathbb{R}}^N$ with finite measure and  $ M_1, M_2>0$. Then there exists $c>0$ depending only on
$\lambda, \mu,     M_1, M_2$ and $|\Omega |$ such that 
estimate (\ref{generalest}) holds
for all $t>0$ and for all diffeomorphisms $\phi$ of class $C^{1,1}$ from $\Omega$ onto an open set $\phi (\Omega )$ in ${\mathbb{R}}^N$ such that inequalities (\ref{fiest1}) are satisfied and  $\delta (\phi )<c^{-1}$.\end{thm}

{\bf Proof.} Let $\phi$ be   diffeomorphism of class $C^{1,1}$ from $\Omega$ onto an open set $\phi (\Omega )$ in ${\mathbb{R}}^N$, satisfying inequalities (\ref{fiest1}). 
Obviously we have
\begin{eqnarray}\label{loclip1}
\lefteqn{    
\left|    
 \frac{Q_{\phi (\Omega)}(C_{\phi }(\beta , w) ) }{\|C_{\phi }(\beta , w) )\|^2_{{\mathcal{L}}^2_t(\phi (\Omega)) }}
-   \frac{Q_{\Omega}(\beta , w ) }{\| (\beta , w) )\|^2_{{\mathcal{L}}^2_t(\Omega ) }}
\right| \le \frac{|Q_{\phi (\Omega)}(C_{\phi }(\beta , w) )- Q_{\Omega}(\beta , w ) |}{ \|C_{\phi }(\beta , w) )\|^2_{{\mathcal{L}}^2_t(\phi (\Omega)) }  }
}\nonumber \\
& &\qquad\qquad\qquad\qquad  +
\frac{  Q_{\Omega }(\beta , w) \left| \| C_{\phi }(\beta , w)   \|^2_{  {\mathcal{L}}^2_t(\phi (\Omega))   }-  \| (\beta , w)   \|^2_{  {\mathcal{L}}^2_t(\Omega )   }\right|      }{ 
 \|C_{\phi }(\beta , w) )\|^2_{{\mathcal{L}}^2_t(\phi (\Omega)) } \|(\beta , w) )\|^2_{{\mathcal{L}}^2_t(\Omega }   }.
\end{eqnarray} 
As in the proof of Lemma~\ref{fiest} , one can prove the existence of a constant $c>0$ depending only on $N, M_1, M_2$ such that 
\begin{equation}\label{loclip2}
\| C_{\phi }(\beta , w)   \|^2_{  {\mathcal{L}}^2_t(\phi (\Omega))   }\geq c  \| (\beta , w)   \|^2_{  {\mathcal{L}}^2_t(\Omega )   } 
\end{equation}
and 
\begin{equation}\label{loclip3}
\left|\| C_{\phi }(\beta , w)   \|^2_{  {\mathcal{L}}^2_t(\phi (\Omega))   }-  \| (\beta , w)   \|^2_{  {\mathcal{L}}^2_t(\Omega )   } \right|\le c \delta(\phi )  \| (\beta , w)   \|^2_{  {\mathcal{L}}^2_t(\Omega )   },
\end{equation}
see also Lemma~\ref{cifi}.   By combining inequalities (\ref{fiest2}) and (\ref{loclip1})-(\ref{loclip3}) we deduce that 
\begin{equation}\label{loclip4}
(1-c\delta (\phi))\frac{Q_{\Omega}(\beta , w ) }{\| (\beta , w) )\|^2_{{\mathcal{L}}^2_t(\Omega ) }}
\le  \frac{Q_{\phi (\Omega)}(C_{\phi }(\beta , w) ) }{\|C_{\phi }(\beta , w) )\|^2_{{\mathcal{L}}^2_t(\phi (\Omega)) }}\le (1+c\delta(\phi ))\frac{Q_{\Omega}(\beta , w ) }{\| (\beta , w) )\|^2_{{\mathcal{L}}^2_t(\Omega ) }}.
\end{equation}
If $1-c\delta (\phi)>0$, it is possible to apply the Min-Max Principle to deduce (\ref{generalest}) from (\ref{loclip4}) combined with Lemma~\ref{cifi}. \hfill $\Box$ \\

\begin{rem}
Since the weak formulation (\ref{reis}) involves only weak derivatives of the first order, one may try  to obtain  stability estimates also under weaker assumptions of $\phi$. For example, one may think of using  bi-Lipschitz domain transformations, i.e., maps  $\phi$ of class $C^{0,1}$ together with their inverses. In this case, one would replace the measure of vicinity $\delta (\phi )$ by the natural weaker measure  of vicinity
$$
\tilde\delta ( \phi)= \| \nabla \phi  -I  \|_{L^{\infty }(\Omega)} .
$$
 In order to prove the corresponding estimate, in the proof of Theorem~\ref{loclip} one should replace the operator $C_{\phi}$ defined in (\ref{cifi0}) by the operator $\tilde C_{\phi }$ defined by 
$$
\tilde C_{\phi }(\beta , w)=(\beta \circ \phi ^{(-1)}, w \circ \phi ^{(-1)}),
$$
for all $(\beta ,w )\in {\mathcal{V}}(\Omega)$.  The definition of the operator $\tilde C_{\phi}$ does not involve $\nabla \phi $ and establishes a linear homeomorphism between ${\mathcal{V}}(\Omega)$ and ${\mathcal{V}}(\phi (\Omega ))$. Unfortunately, the summand $ \int_{\Omega}(\nabla w -\beta )\cdot (\nabla v-\eta)dx$ in the quadratic form (\ref{reis}) does not behave well under the transformation  $\tilde C_{\phi}$ and this would lead to an estimate depending  on $t$. Namely, one would obtain the estimate
\begin{equation}\label{generalestrem}
|\gamma_{n,t}[\phi (\Omega)]-\gamma_{n,t}[\Omega]|\le \frac{c}{t^2} \gamma_{n,t}[\Omega] \tilde \delta (\phi ),
\end{equation}
where the presence of a better measure of vicinity $\tilde\delta (\phi) $ is compensated by the presence of the factor $t^2$ which spoils the estimate for $t$ close to zero. 

In any case, using domain transformations $\phi$ of class $C^{1,1}$ and the corresponding strong measure of vicinity $\delta (\phi)$ is enough for our purpose of obtaining estimates via Hausdorff distance. 
 \end{rem}

\subsection{Estimates via atlas and Hausdorff distance}

In general, even if  two open sets $\Omega_1$ and $\Omega_2$ are known to be diffeomorphic, it is not easy to construct a diffeomorphism $\phi$ such that $\phi (\Omega_1)=\Omega_2$  and provide information on  $\delta (\phi)$ in terms of  explicit geometric quantities. However, if $\Omega_1$, $\Omega_2$
 belong to the same class $C({\mathcal{A}})$ then it is possible to construct a suitable diffeomorphism $\phi $ such that $\phi (\Omega_1)\subset \Omega_2$ and estimate  $\delta (\phi)$  via the atlas distance (\ref{dev1}).  Such construction was first used in Burenkov and Davies~\cite{buda} and then implemented in \cite{bulahigh}. We briefly recall it. 

Let ${\mathcal{A}}  =(\rho ,s,s', \left\{V_j\right\}_{j=1}^s,
\left\{{r}_j\right\}_{j=1}^s)$ be an atlas in ${\mathbb{R}}^N$
and let  $\{\psi_j\}_{j=1}^{s}$  be a partition of unity such that $\psi_j\in C^{\infty}_c(\mathbb{R}^N)$, ${\rm supp}\,\psi_j\subset {(V_j)}_{\frac{3}{4}\rho}$, $0\le \psi_j \le 1$ and  $\sum_{j=1}^s\psi_j(x)=1$ for all $x\in \cup_{j=1}^s(V_j)_{\rho}$. 
For $\epsilon \geq 0$ we consider  the following transformation
\begin{equation}
\label{budaout1}
\phi_{\epsilon}(x)=x-\epsilon\sum_{j=1}^s\xi_j\psi_j(x)\, ,\ \ \ x\in\mathbb{R}^N,
\end{equation}
where
$\xi_j={r}_j^{(-1)}((0,\dots ,1))$.

Then we recall the following technical lemma from  \cite{bulahigh}.
\begin{lemm}
\label{budaout}
Let ${\mathcal{A}}$ be an atlas in ${\mathbb{R}}^N$.
Then there exist $M, M_1, M_2, E>0$ depending only on $N$ and ${\mathcal{A}}  $
 such that   $\phi_{\epsilon }$  satisfies   (\ref{fiest1}) and such that $\delta (\phi_{\epsilon})\le M\epsilon $   for all $\epsilon \in [0, E[$.
 Moreover,  $\phi_{\epsilon }(\Omega_1)\subset \Omega_2$
for all $\epsilon \in [0,E[ $ and 
for all $\Omega_1,\Omega_2\in C({\mathcal{A}}) $ such that $\Omega_2\subset \Omega_1$ and
$d_{{\mathcal{A}}}(\Omega_1,\Omega_2)<\epsilon /s$.
\end{lemm}

Proceeding as in \cite{bulahigh} we can prove the following

\begin{thm}\label{dirthm}
Let ${\mathcal{A}}$ be an atlas in ${\mathbb{R}}^N$. 
 Then  there exists $c >0$ depending only on
$  {\mathcal{A}}, \lambda, \mu$   such that estimate (\ref{dirthm1}) holds
 for all $n\in {\mathbb{N}}$, $t>0$ and for all $\Omega_1, \Omega_2\in C({\mathcal{A}})$ satisfying $d_{\mathcal{A}}(\Omega_1,\Omega_2) <c^{-1} $.
\end{thm}

{\bf Proof.} Let $E>0$ be as in Lemma~\ref{budaout} and let $\Omega_1,\Omega_2\in C({\mathcal{A}}) $ be  such that
$d_{{\mathcal{A}}}(\Omega_1,\Omega_2)<\epsilon /s$. Clearly $\Omega_1\cap\Omega_2\in C({\mathcal{A}}) $ and $d_{{\mathcal{A}}}(\Omega_1\cap\Omega_2,\Omega_1), d_{{\mathcal{A}}}(\Omega_1\cap\Omega_2,\Omega_1)<\epsilon /s$. Thus  by Lemma~\ref{budaout} we have that $\phi_{\epsilon }(\Omega_1), \phi_{\epsilon }(\Omega_2)\subset \Omega_1\cap\Omega_2$. By the monotonicity of the eigenvalues with respect to inclusion, we immediately get
\begin{equation}
\gamma_{n,t}[\Omega_i]\le \gamma_{n,t}[\Omega_1\cap \Omega_2]\le \gamma_{n,t}[\phi_{\epsilon}(\Omega _i)],
\end{equation}
for all $i=1,2$. Moreover, by combining Theorem~\ref{loclip} and Lemma~\ref{budaout}, we deduce that there exists $c$ as in the statement such that 
\begin{equation}\label{dirthmlast} | \gamma_{n,t}[\Omega_i]- \gamma_{n,t}[\Omega_1\cap \Omega_2]|\le 
| \gamma_{n,t}[\phi_{\epsilon }(\Omega_i)] -\gamma_{n,t}[\Omega_i]  |\le c \gamma_{n,t}[\Omega_i] \epsilon,
\end{equation}
for all $i=1,2$, provided $\epsilon \le c^{-1}$. Inequality (\ref{dirthm1}) easily follows by choosing $\epsilon =2sd_{{\mathcal{A}}}(\Omega_1,\Omega_2)$ in (\ref{dirthmlast}). \hfill $\Box$. \\

We note that by Theorem~\ref{metricsurv} and estimate (\ref{dirthm1}), it immediately follows that if $\omega $ is a modulus of continuity  as in Definition~\ref{atclass} then there exist $c >0$ depending only on
$  {\mathcal{A}}, \omega, \lambda, \mu $   such that
\begin{equation}
\label{dirthm2}
| \gamma_{n,t}[\Omega_1]-\gamma_{n  ,t}[\Omega_2] |\le  c\max\{ \gamma_{n,t}[\Omega_1] , \gamma_{n  ,t}[\Omega_2]\}  \omega (d_{{\mathcal{H}}   {\mathcal{P}}}(\partial \Omega_1,\partial \Omega_2)),
\end{equation}
 for all $n\in {\mathbb{N}}$, $t>0$ and  for all  $\Omega_1, \Omega_2 \in C^{\omega (\cdot ) }_M({\mathcal{A}})$ satisfying the condition $d_{{\mathcal{H}}   {\mathcal{P}}}(\Omega_1,\Omega_2) <c^{-1} $.\\

In several papers devoted to stability estimates for domain perturbation problems, the vicinity of two domains is  described by means of $\epsilon$-neigh\-bor\-hoods of the boundaries defined by
the Euclidean distance, see e.g.,  \cite{buda} and Davies~\cite{daviespaper}.  This can be done also in the case of the Reissner-Mindlin system. Indeed, one can  prove the following

\begin{cor}\label{inclestthm}
Let ${\mathcal{A}}$ be an atlas in ${\mathbb{R}}^N$, $\omega$ a modulus of continuity as in Definition~\ref{atclass} and $ M>0$. 
 Then  there exists $c >0$ depending only on
$   {\mathcal{A}}, \omega, \lambda, \mu,   M $   such that
\begin{equation}\label{inclest}
| \gamma_{n,t}[\Omega_1]-\gamma_{n  ,t}[\Omega_2] |\le  c\max\{ \gamma_{n,t}[\Omega_1] , \gamma_{n  ,t}[\Omega_2]\}  \omega (\epsilon),
\end{equation}
for all $n\in {\mathbb{N}}$, $t>0$, $\epsilon \in ]0, c^{-1}[$ and for all $\Omega_1, \Omega_2\in C^{\omega (\cdot ) }_M({\mathcal{A}})$ such that 
\begin{equation}\label{incl}
(\Omega_1)_{\epsilon }\subset \Omega_2\subset (\Omega_1)^{\epsilon},\ \ {\rm or}\ \ (\Omega_2)_{\epsilon }\subset \Omega_1\subset (\Omega_2)^{\epsilon}.
\end{equation}
\end{cor}

{\bf Proof. } Note that if $\Omega_1$ and $\Omega_2$ satisfy one of the inclusions in (\ref{incl}) then  $d_{{\mathcal{H}}   {\mathcal{P}}}(\partial \Omega_1,\partial \Omega_2)\le \epsilon$, which combined with inequality (\ref{dirthm2}) allows to deduce (\ref{inclest}). \hfill $\Box$\\

\section{Shape differentiability}

Given a bounded open set in ${\mathbb{R}}^N$, we denote by $C^{0,1}(\Omega ; {\mathbb{R}}^N )$ the set of Lipschitz continuous maps from $\Omega $ to ${\mathbb{R}}^N$.
By  ${\rm BLip}(\Omega)$ we denote the set of functions $\phi \in C^{0,1}(\Omega ; {\mathbb{R}}^N )$ such that $\phi $ is injective and the inverse $\phi ^{(-1)}:\phi (\Omega)\to \Omega$ is  Lipschitz continuous. We shall think of $C^{0,1}(\Omega ; {\mathbb{R}}^N )$ as a Banach space endowed with the standard norm defined by
$$
\| \phi \|_{C^{0,1}(\Omega ;{\mathbb{R}}^N)} = \| \phi \|_{L^{\infty }(\Omega )}+{\rm Lip}(\phi),
$$
for all $\phi \in C^{0,1}(\Omega ;{\mathbb{R}}^N)$, where ${\rm Lip}(\phi)$ is the Lipschitz constant of $\phi$.  We recall   that ${\rm BLip}(\Omega)$ is an open set in $C^{0,1}(\Omega ; {\mathbb{R}}^N )$, see e.g., \cite[Lemma~3.11]{lala}.

In this section, we prove analyticity results for the maps  $\phi \mapsto \gamma_{n,t}[\phi (\Omega )]$, defined for $\phi\in {\rm BLip}(\Omega)$.  
To shorten our notation, in the sequel we shall write $\gamma_{n,t}[\phi] $ instead of $\gamma_{n, t}[\phi(\Omega)]$.

As is known, when dealing with differentiability properties of the eigenvalues, it is necessary to pay attention to bifurcation phenomena associated with multiple eigenvalues.
 Following \cite{lala, lalcri}, 
given a finite non-empty subset of $\mathbb{N}$, we set
	$$\mathcal{A}_{F,t}(\Omega)=\{\phi\in {\rm BLip}(\Omega):\gamma_{l,t}[\phi]\notin\{\gamma_{j,t}[\phi]:j\in F\}\ \forall l\in\mathbb{N}
				\setminus F \}$$
and 
$$\Theta_{F,t}(\Omega)=\{\phi\in\mathcal{A}_{F,t}(\Omega):\gamma_{j,t}[\phi]\ \text{have a common value}\ \gamma_{F,t}[\phi]\ \forall j\in F\}.$$
	
Then we can prove the following real-analyticity result in the spirit of the results in \cite{buoso, buosoplates, lala}.  

	\begin{thm}
	\label{thesame}
	Let $\Omega$ be a bounded open set in $\mathbb{R}^N$, $t>0$ and $F$ a finite non-empty subset of $\mathbb{N}$.
	Then the following statements hold.
	\begin{enumerate}[i)]
	\item The set $\mathcal{A}_{F,t}(\Omega)$ is open in $C^{0,1} (\Omega ; {\mathbb{R}}^N)$. Moreover, for every $s\in\{1,\dots,|F|\}$ the real valued function $\Gamma_{F,t}^{(s)}$ defined on $\mathcal{A}_{F,t}(\Omega)$ by
		$$\Gamma_{F,t}^{(s)}[\phi]=\sum_{\substack{j_1<\dots<j_s\\ j_1,\dots , j_s\in F}}\gamma_{j_1,t}[\phi]\cdots\gamma_{j_s,t}[\phi]$$
for all $\phi \in \mathcal{A}_{F,t}(\Omega)$, is real analytic.
\item   Let $\tilde{\phi}\in\Theta_{F,t}(\Omega)$ be such that $\tilde{\phi}(\Omega)$ is of class 
		$ C^{1,1}$. Then for every $s\in\{1,\dots,|F|\}$ the Fr\'{e}chet differential of the function $\Gamma_{F,t}^{(s)}$ at the point $\tilde \phi$ is provided by 
the formula 
		\begin{multline}
		\label{deriv}
			d_{|_{\phi=\tilde{\phi}}}\Gamma_{F,t}^{(s)}[\psi]=-
			\gamma_{F,t}^{s-1}[\tilde{\phi}]\binom{|F|-1}{s-1}
			\sum_{l=1}^{|F|}
			\int_{\partial\tilde{\phi}(\Omega)}\left(
		\frac{\mu}{12}\left|\frac{\partial\beta^{(l)}}{\partial n}\right|^2\right.\\
	\left.+\frac{\mu+\lambda}{12}\left(\frac{\partial\beta^{(l)}}{\partial n}\cdot n\right)^2
	+\frac{\mu k}{t^2}\left(\frac{\partial w^{(l)}}{\partial n}\right)^2\right)\zeta\cdot n d\sigma,
		\end{multline}
		for all $\psi\in C^{0,1}(\Omega ; {\mathbb{R}}^N)$, where $\zeta=\psi\circ\tilde{\phi}^{(-1)}$ and  $(\beta^{(1)}, w^{(1)}),$ $ \dots,$ $ (\beta^{(|F|)},$ $ w^{(|F|)})$
		is an orthonormal basis in ${\mathcal{L}}^2_t(\tilde\phi (\Omega ))$ for the eigenspace associated with $\gamma_{F,t}[\tilde{\phi}]$.
	\end{enumerate} 
	\end{thm}

{\bf Proof.} The proof can be deduced by the abstract results in \cite{lala} as follows.  We consider the operator $R_{\phi (\Omega),t}$ as an operator acting from the space ${\mathcal{V}}(\phi (\Omega))$ to its dual and we pull-it back to $\Omega$ by
changing variables via $\phi$. Namely, the pull-back ${\mathcal{R}}_{\phi ,t }$ of $R_{\phi (\Omega),t}$ is the operator  defined from  ${\mathcal{V}}(\Omega)$ to its dual which takes any  
 $(\theta ,u)\in {\mathcal{V}}(\Omega)$ to the functional ${\mathcal{R}}_{\phi ,t}(\theta ,u)$ defined by 
\begin{eqnarray}\label{pul} \lefteqn{
	{\mathcal{R}}_{\phi ,t}(\theta , u)( \dot \theta , \dot u)=
		\frac{\mu}{12}\int_{\Omega}\left(\nabla(\theta\circ\phi^{(-1)}):\nabla(\dot \theta\circ\phi^{(-1)})\right)\circ\phi|\det D\phi|dx}\nonumber \\
	& &\quad 	+\frac{\mu+\lambda}{12}\int_{\Omega}\left(\mathrm{div}(\theta\circ\phi^{(-1)})\mathrm{div}(\dot\theta\circ\phi^{(-1)})\right)\circ\phi|\det D\phi|dx\nonumber \\
	& &\quad   +\frac{\mu k}{t^2}\int_{\Omega}(\nabla (u\circ\phi^{(-1)})\circ\phi-\theta)\cdot(\nabla (\dot u\circ\phi^{(-1)})\circ\phi-\dot \theta)|\det D\phi|dx,
	\end{eqnarray}
	for all $( \dot \theta , \dot u )\in {\mathcal{V}}(\Omega)$. Similarly, we consider the map $\mathcal{J}_{\phi ,t}$ from ${\mathcal{V}}(\Omega)$ to its dual  defined by
	\begin{equation}\label{geifi}\mathcal{J}_{\phi ,t}(\theta ,u)( \dot\theta , \dot u)=\int_{\Omega}\left(u\dot u+\frac{t^2}{12}\theta\cdot\dot \theta\right)|\det D\phi|dx,\end{equation}
	for all $(\theta , u),(\dot\theta , \dot u)\in {\mathcal{V}}(\Omega)$. 
	Note that ${\mathcal{R}}_{\phi ,t}(\theta , u)(\dot \theta , \dot u)$ can be considered as a  scalar product in ${\mathcal{V}}(\Omega)$ and the corresponding norm is equivalent to the standard Sobolev norm.  Accordingly, we can think of ${\mathcal{V}}(\Omega)$  as a Hilbert space endowed with such scalar product. Thus, by the Riesz Representation Theorem applied to  ${\mathcal{V}}(\Omega)$, it follows that the operator ${\mathcal{R}}_{\phi ,t}$ is invertible.

 It is easy to see that $( \beta ,w )\in {\mathcal{V}}(\phi (\Omega))$ is an eigenvector associated with an eigenvalue  $\gamma$  of the operator $R_{\phi (\Omega), t}$ if and only if $ {\mathcal{R}}_{\phi ,t}(\beta\circ \phi , w\circ \phi )=\gamma \mathcal{J}_{\phi ,t}(\beta\circ \phi , w\circ \phi )$. This implies that the eigenvalues of the operator  $R_{\phi (\Omega),t}$ are the reciprocal of the eigenvalues of the operator $T_{\phi ,t}$ defined from ${\mathcal{V}}(\Omega)$ to itself by 
\begin{equation}\label{pul1}T_{\phi ,t}={\mathcal{R}}^{(-1)}_{\phi ,t}\circ \mathcal{J}_{\phi ,t}. \end{equation}
It turns out that $T_{\phi ,t}$ is a compact self-adjoint operator on the Hilbert space ${\mathcal{V}}(\Omega )$.  Note that the  operators ${\mathcal{R}}_{\phi ,t}$, $\mathcal{J}_{\phi ,t }$, hence $T_{\phi ,t}$ depend real-analytically on $\phi$, since they are obtained as composition of real-analytic maps. Thus, it is possible to apply the general results in \cite{lala} and conclude that the elementary symmetric functions $\sum_{j_1<\dots<j_s\in F}\gamma^{-1}_{j_1,t}[\phi]\cdots\gamma^{-1}_{j_s,t}[\phi]$  of the eigenvalues of $T_{\phi ,t}$ depend real-analytically on $\phi$. Then, by arguing as in \cite{lala}, one can easily deduce the validity of statement (i). 

As for statement (ii),  	we set  $\theta^{(l)}=\beta^{(l)}\circ\tilde{\phi}$ and $u^{(l)}=w^{(l)}\circ\tilde{\phi}$ for $l=1,\dots,|F|$.   By arguing as in \cite{lala} we obtain

	\begin{equation}\nonumber
	d|_{\phi=\tilde{\phi}}(\Gamma_{F,t}^{(s)})[\psi]
	=-\gamma_{F,t}^{s}[\tilde{\phi}]\binom{|F|-1}{s-1}\sum_{l=1}^{|F|}{\mathcal{R}}_{\tilde{\phi}, t}\left(
		d|_{\phi=\tilde{\phi}}T_{t,\phi}[\psi](\theta^{(l)}, u^{(l)})\right)  \left( (\theta^{(l)}, u^{(l)})\right).
	\end{equation}
	Then one can easily prove formula (\ref{deriv}) using Lemma \ref{tomare} below. \hfill $\Box$ \\

	\begin{lemm}
	\label{tomare}
	Let $\Omega$ be a bounded open set in $\mathbb{R}^N$ and $\tilde \phi \in  {\rm BLip}(\Omega)$ be such that $\tilde \phi (\Omega)$ is of class $C^{1,1}$. Let $t>0$ and  $( \beta^{(i)}, w^{(i)})\in {\mathcal{V}}(\tilde \phi (\Omega ))$, $i=1,2$ be eigenvectors associated with an eigenvalue $\tilde\gamma$ of the operator $R_{\tilde\phi (\Omega),t}$. Let $\theta^{(i)}=\beta^{(i)}\circ\tilde \phi$, $u^{(i)}=w^{(i)}\circ \tilde\phi$,  $i=1,2$.
Then we have
	\begin{multline}
	{\mathcal R}_{\tilde{\phi},t}\left(d|_{\phi=\tilde{\phi}}T_{\phi ,t}[\psi]
	(\theta^{(1)} , u^{(1)})\right)(\theta^{(2)}, u^{(2)})	=\tilde\gamma^{-1}\frac{\mu}{12}\int_{\partial\tilde{\phi}(\Omega)}\frac{\partial\beta^{(1)}}{\partial n}\cdot\frac{\partial\beta^{(2)}}{\partial n}\zeta\cdot n d\sigma \\ \qquad
	+\tilde\gamma^{-1}\frac{\mu+\lambda}{12}\int_{\partial\tilde{\phi}(\Omega)}\frac{\partial\beta^{(1)}}{\partial n}\cdot n\frac{\partial\beta^{(2)}}{\partial n}\cdot n\zeta\cdot n d\sigma
	+\tilde\gamma^{-1}\frac{\mu k}{t^2}\int_{\partial\tilde{\phi}(\Omega)}\frac{\partial w^{(1)}}{\partial n}\frac{\partial w^{(2)}}{\partial n}\zeta\cdot n d\sigma,
	\end{multline}
	for all $\psi\in C^{0,1}(\Omega; {\mathbb{R}}^N)$, where $\zeta=\psi\circ\tilde{\phi}^{(-1)}$ and ${\mathcal R}_{\tilde{\phi},t}$, $T_{\phi ,t}$ are defined by (\ref{pul}), (\ref{pul1}) respectively.
	\end{lemm}

{\bf Proof.} First of all, we note that  by classical regularity theory, the eigenvectors $(\beta^{(i)}, w^{(i)})$, $i=1,2$, belong to
	$(H^2(\tilde{\phi}(\Omega)))^{N}\times H^2(\tilde{\phi}(\Omega))$. This will be used in most of the following computations.

	By standard calculus in normed spaces we have
	\begin{multline}\label{maybefirst}
			{\mathcal{R}}_{\tilde{\phi},t}\left[\mathrm{d}|_{\phi=\tilde{\phi}}\left(
				{\mathcal{R}}_{\phi ,t}^{(-1)}\circ\mathcal{J}_{\phi ,t}[\psi](\theta^{(1)}, u^{(1)}), (\theta^{(2)}, u^{(2)})\right)\right]  \\ 
			= {\mathcal{R}}_{\tilde{\phi},t}\left[{\mathcal{R}}_{\tilde{\phi},t}^{(-1)}\circ
							\mathrm{d}|_{\phi=\tilde{\phi}}\mathcal{J}_{\phi ,t}[\psi](\theta^{(1)}, u^{(1)}), (\theta^{(2)}, u^{(2)})\right]\\
			 						+{\mathcal{R}}_{\tilde{\phi},t}\left[\mathrm{d}|_{\phi=\tilde{\phi}}
				{\mathcal{R}}_{\phi,t}^{(-1)}[\psi]\circ\mathcal{J}_{\tilde{\phi} ,t}(\theta^{(1)}, u^{(1)}), (\theta^{(2)}, u^{(2)})\right]	.
	\end{multline}

Now we note that
	\begin{multline}
	\label{second}
		{\mathcal{R}}_{\tilde{\phi},t}\left[\mathrm{d}|_{\phi=\tilde{\phi}}
				{\mathcal{R}}_{\phi ,t}^{(-1)}[\psi]\circ\mathcal{J}_{\tilde{\phi},t}(\theta^{(1)}, u^{(1)}), (\theta^{(2)}, u^{(2)})\right]\\
		=-{\mathcal{R}}_{\tilde{\phi},t}\left[{\mathcal{R}}_{\tilde{\phi},t}^{(-1)}\circ\mathrm{d}|_{\phi=\tilde{\phi}}{\mathcal{R}}_{\phi,t}[\psi]\circ {\mathcal{R}}_{\tilde{\phi},t}^{(-1)}\circ\mathcal{J}_{\tilde{\phi},t}(\theta^{(1)}, u^{(1)}), (\theta^{(2)}, u^{(2)})\right]\\
				=-\tilde{\gamma}^{-1}\left(d|_{\phi=\tilde{\phi}}{\mathcal{R}}_{\phi,t}[\psi](\theta^{(1)}, u^{(1)})\right)(\theta^{(2)}, u^{(2)} ).
	\end{multline}
	
	By standard calculus we have
	\begin{equation}
	\label{determinante}
	\left[\left(d|_{\phi=\tilde{\phi}}(\det \nabla\phi)[\psi]\right)\circ\tilde{\phi}^{(-1)}\right]\det \nabla \tilde{\phi}^{(-1)}=
		\mathrm{div}\zeta,
	\end{equation}
	hence
	\begin{equation}\label{degei}(d|_{\phi=\tilde{\phi}}\mathcal{J}_{\phi ,t}[\psi][(\theta^{(1)} , u^{(1)}  )])[(\theta^{(2)}, u^{(2)}  )]=
		\int_{\tilde{\phi}(\Omega)}\left(w^{(1)}w^{(2)}+\frac{t^2}{12}\beta^{(1)}\beta^{(2)}\right)\mathrm{div}\zeta dy.\end{equation}
	
Note that, in order to shorten our notation, in the sequel summation symbols will be omitted. By standard calculus in normed space and  changing variables we get
	\begin{multline}
	\label{important}
		\left(d|_{\phi=\tilde{\phi}}{\mathcal{R}}_{t,\phi}[\psi](\theta^{(1)}, u^{(1)})\right)(\theta^{(2)}, u^{(2)})\\
		=-\frac{\mu}{12}\int_{\tilde{\phi}(\Omega)}\left(\frac{\partial\beta^{(1)}_i}{\partial y_r}\frac{\partial\beta^{(2)}_i}{\partial y_j}
				+\frac{\partial\beta^{(2)}_i}{\partial y_r}\frac{\partial\beta^{(1)}_i}{\partial y_j}\right)
					\frac{\partial\zeta_r}{\partial y_j}dy
		+\frac{\mu}{12}\int_{\tilde{\phi}(\Omega)}\frac{\partial\beta^{(1)}_i}{\partial y_j}\frac{\partial\beta^{(2)}_i}{\partial y_j}\mathrm{div}\zeta dy\\
		-\frac{\mu+\lambda}{12}\int_{\tilde{\phi}(\Omega)}\left(\frac{\partial\beta^{(1)}_i}{\partial y_r}\mathrm{div}\beta^{(2)}+
			\frac{\partial\beta^{(2)}_i}{\partial y_r}\mathrm{div}\beta^{(1)}\right)\frac{\partial\zeta_r}{\partial y_i}dy\\
		+\frac{\mu+\lambda}{12}\int_{\tilde{\phi}(\Omega)}\mathrm{div}\beta^{(1)}\mathrm{div}\beta^{(2)}\mathrm{div}\zeta dy
		-\frac{\mu k}{t^2}\int_{\tilde{\phi}(\Omega)}\frac{\partial w^{(1)}}{\partial y_r}\frac{\partial\zeta_r}{\partial y_i}
			\left(\frac{\partial w^{(2)}}{\partial y_i}-\beta^{(2)}_i\right)dy\\
		-\frac{\mu k}{t^2}\int_{\tilde{\phi}(\Omega)}\left(\frac{\partial w^{(1)}}{\partial y_i}-\beta^{(1)}_i\right)
		\frac{\partial w^{(2)}}{\partial y_r}\frac{\partial\zeta_r}{\partial y_i}dy\\
		+\frac{\mu k}{t^2}\int_{\tilde{\phi}(\Omega)}\left(\frac{\partial w^{(1)}}{\partial y_i}-\beta^{(1)}_i\right)
		\left(\frac{\partial w^{(2)}}{\partial y_i}-\beta^{(2)}_i\right)\mathrm{div}\zeta dy.			
	\end{multline}

	Now note that
	\begin{multline}\label{importantbis}
	\int_{\tilde{\phi}(\Omega)}\frac{\partial\beta^{(1)}_i}{\partial y_r}\frac{\partial\beta^{(2)}_i}{\partial y_j}\frac{\partial\zeta_r}{\partial y_j}dy
	=\int_{\partial\tilde{\phi}(\Omega)}\frac{\partial\beta^{(1)}_i}{\partial n}\frac{\partial\beta^{(2)}_i}{\partial n}\zeta\cdot nd\sigma\\
	-\int_{\tilde{\phi}(\Omega)}\Delta\beta^{(2)}\cdot (\nabla\beta^{(1)}\cdot\zeta )dy
	-\int_{\tilde{\phi}(\Omega)}\frac{\partial\beta^{(2)}_i}{\partial y_j}\frac{\partial^2\beta^{(1)}_i}{\partial y_j\partial y_r}\zeta_r dy\\
	=-\int_{\tilde{\phi}(\Omega)}\Delta\beta^{(2)}\cdot (\nabla\beta^{(1)}\cdot\zeta )dy
	+\int_{\tilde{\phi}(\Omega)}\frac{\partial\beta^{(1)}_i}{\partial y_j}\frac{\partial^2\beta^{(2)}_i}{\partial y_j\partial y_r}\zeta_r dy\\
	+\int_{\tilde{\phi}(\Omega)}\frac{\partial\beta^{(1)}_i}{\partial y_j}\frac{\partial\beta^{(2)}_i}{\partial y_j}\mathrm{div}\zeta dy.	
	\end{multline}

Note that here and in the sequel we also use the fact that if $U$ is a smooth open set and  $f\in H^2(U)\cap H^1_0(U)$ then $\nabla f=\frac{\partial f}{\partial n}n$ on $\partial U$; moreover,   if $g\in (H^2(U)\cap H^1_0(U))^N $ then ${\rm div} g=\frac{\partial g}{\partial n}\cdot n$ on $\partial U$.   

	By (\ref{importantbis}) the sum of the first two integrals in  the right-hand side of (\ref{important}) equals
	\begin{multline}
	\label{primi3}
	-\frac{\mu}{12}\int_{\partial\tilde{\phi}(\Omega)}\frac{\partial\beta^{(1)}}{\partial n}\cdot\frac{\partial\beta^{(2)}}{\partial n}\zeta\cdot nd\sigma\\
	+\frac{\mu}{12}\int_{\tilde{\phi}(\Omega)}\left(\Delta\beta^{(1)}_i\nabla\beta^{(2)}_i+\Delta\beta^{(2)}_i\nabla\beta^{(1)}_i\right)\cdot\zeta dy.
	\end{multline}	
		
	Now we observe that
	\begin{multline}
	\int_{\tilde{\phi}(\Omega)}\frac{\partial\beta^{(1)}_i}{\partial y_r}\frac{\partial\zeta_r}{\partial y_i}\mathrm{div}\beta^{(2)}dy
	=\int_{\partial\tilde{\phi}(\Omega)}\frac{\partial\beta^{(1)}}{\partial n}\cdot n\mathrm{div}\beta^{(2)}\zeta\cdot n d\sigma\\
	-\int_{\tilde{\phi}(\Omega)}\frac{\partial\mathrm{div}\beta^{(1)}}{\partial y_r}\zeta_r\mathrm{div}\beta^{(2)}dy
	-\int_{\tilde{\phi}(\Omega)}\frac{\partial\mathrm{div}\beta^{(2)}}{\partial y_i}\frac{\partial\beta^{(1)}_i}{\partial y_r}\zeta_rdy\\
	=
	-\int_{\tilde{\phi}(\Omega)}\frac{\partial\mathrm{div}\beta^{(2)}}{\partial y_i}\frac{\partial\beta^{(1)}_i}{\partial y_r}\zeta_rdy
	 +\int_{\tilde{\phi}(\Omega)}\mathrm{div}\beta^{(1)}\mathrm{div}\beta^{(2)}\mathrm{div}\zeta dy\\
	+\int_{\tilde{\phi}(\Omega)}\mathrm{div}\beta^{(1)}\frac{\partial\mathrm{div}\beta^{(2)}}{\partial y_r}\zeta_rdy.
	\end{multline}
	
Thus,  the sum of third and the fourth integral  in the right-hand side of (\ref{important}) is equal to
	\begin{multline}
	\label{secondi3}
	\frac{\mu+\lambda}{12}\int_{\tilde{\phi}(\Omega)}\left(\frac{\partial\mathrm{div}\beta^{(1)}}{\partial y_i}\frac{\partial\beta^{(2)}_i}{\partial y_r}
	+\frac{\partial\mathrm{div}\beta^{(2)}}{\partial y_i}\frac{\partial\beta^{(1)}_i}{\partial y_r}\right)\zeta_rdy\\
	-\frac{\mu+\lambda}{12}\int_{\partial\tilde{\phi}(\Omega)}\frac{\partial\beta^{(1)}}{\partial n}\cdot n\frac{\partial\beta^{(2)}}{\partial n}\cdot n
		\zeta\cdot n d\sigma.
	\end{multline}
	
	Now note that
	\begin{multline}\label{swap}
	\int_{\tilde{\phi}(\Omega)}\frac{\partial w^{(1)}}{\partial y_r}\frac{\partial\zeta_r}{\partial y_i}\left(\frac{\partial w^{(2)}}{\partial y_i}-\beta^{(2)}_i\right)dy
	=\int_{\partial\tilde{\phi}(\Omega)}\frac{\partial w^{(1)}}{\partial n}\frac{\partial w^{(2)}}{\partial n}\zeta\cdot n d\sigma\\
	-\int_{\tilde{\phi}(\Omega)}\frac{\partial w^{(1)}}{\partial y_r}\zeta_r\left(\Delta w^{(2)}-\mathrm{div}\beta^{(2)}\right)dy
	-\int_{\tilde{\phi}(\Omega)}\frac{\partial^2w^{(1)}}{\partial y_i\partial y_r}\zeta_r\left(\frac{\partial w^{(2)}}{\partial y_i}-\beta^{(2)}_i\right)dy\\
	=-\int_{\tilde{\phi}(\Omega)}\frac{\partial w^{(1)}}{\partial y_r}\zeta_r\left(\Delta w^{(2)}-\mathrm{div}\beta^{(2)}\right)dy
	+\int_{\tilde{\phi}(\Omega)}\nabla w^{(1)}(\nabla w^{(2)}-\beta^{(2)})\mathrm{div}\zeta dy\\
	+\int_{\tilde{\phi}(\Omega)}\frac{\partial w^{(1)}}{\partial y_i}
		\left(\frac{\partial^2w^{(2)}}{\partial y_i\partial y_r}-\frac{\partial\beta^{(2)}_i}{\partial y_r}\right)\zeta_r dy.
	\end{multline}
	
	By using the second equality in (\ref{swap}), and  the first equality in (\ref{swap}) with $(\beta ^{(1)}, w^{(1)})$ replaced by $(\beta ^{(2)}, w^{(2)})$, we get that  the sum of the last three integrals  in (\ref{important}) is equal to
	\begin{multline}
	\label{ultimi3}
	-\frac{\mu k}{t^2}\int_{\partial\tilde{\phi}(\Omega)}\frac{\partial w^{(1)}}{\partial n}\frac{\partial w^{(2)}}{\partial n}\zeta\cdot n d\sigma\\
	+\frac{\mu k}{t^2}\int_{\tilde{\phi}(\Omega)}(\Delta w^{(1)}-\mathrm{div}\beta^{(1)})\nabla w^{(2)}\cdot\zeta dy
	+\frac{\mu k}{t^2}\int_{\tilde{\phi}(\Omega)}(\Delta w^{(2)}-\mathrm{div}\beta^{(2)})\nabla w^{(1)}\cdot\zeta dy\\
	-\frac{\mu k}{t^2}\int_{\tilde{\phi}(\Omega)}\beta^{(1)}(\nabla w^{(2)}-\beta^{(2)})\mathrm{div}\zeta dy
	+\frac{\mu k}{t^2}\int_{\tilde{\phi}(\Omega)}\frac{\partial w^{(1)}}{\partial y_i}\frac{\partial\beta^{(2)}_i}{\partial y_r}\zeta_rdy\\
	-\frac{\mu k}{t^2}\int_{\tilde{\phi}(\Omega)}\beta^{(1)}_i\frac{\partial^2w^{(2)}}{\partial y_i\partial y_r}\zeta_rdy
	=-\frac{\mu k}{t^2}\int_{\partial\tilde{\phi}(\Omega)}\frac{\partial w^{(1)}}{\partial n}\frac{\partial w^{(2)}}{\partial n}\zeta\cdot n d\sigma\\
	+\frac{\mu k}{t^2}\int_{\tilde{\phi}(\Omega)}(\Delta w^{(1)}-\mathrm{div}\beta^{(1)})\nabla w^{(2)}\cdot\zeta dy
	+\frac{\mu k}{t^2}\int_{\tilde{\phi}(\Omega)}(\Delta w^{(2)}-\mathrm{div}\beta^{(2)})\nabla w^{(1)}\cdot\zeta dy\\
	+\frac{\mu k}{t^2}\int_{\tilde{\phi}(\Omega)}\left(\frac{\partial w^{(1)}}{\partial y_i}-\beta^{(1)}_i\right)\frac{\partial\beta^{(2)}_i}{\partial y_r}\zeta_r dy
	+\frac{\mu k}{t^2}\int_{\tilde{\phi}(\Omega)}\left(\frac{\partial w^{(2)}}{\partial y_i}-\beta^{(2)}_i\right)\frac{\partial\beta^{(1)}_i}{\partial y_r}\zeta_r dy.
	\end{multline}
	
	Using the fact that
	$$-\frac{\mu}{12}\Delta\beta^{(i)}-\frac{\mu+\lambda}{12}\nabla\mathrm{div}\beta^{(i)}  -\frac{\mu k}{t^2}(\nabla w^{(i)}-\beta^{(i)})
		=\frac{\tilde{\gamma}t^2}{12}\beta^{(i)},$$	
and
$$-\frac{\mu k}{t^2}(\Delta w^{(i)}-  \mathrm{div}\beta^{(i)})=\tilde{\gamma} w^{(i)},$$
 for $i=1,2$, we get that
	\begin{multline}
	\left(d|_{\phi=\tilde{\phi}}{\mathcal{R}}_{t,\phi}[\psi](\theta^{(1)}, u^{(1)})\right)(\theta^{(2)}, u^{(2)})\\
	=-\frac{\mu}{12}\int_{\partial\tilde{\phi}(\Omega)}\frac{\partial\beta^{(1)}}{\partial n}\cdot\frac{\partial\beta^{(2)}}{\partial n}\zeta\cdot n d\sigma
	-\frac{\mu+\lambda}{12}\int_{\partial\tilde{\phi}(\Omega)}\frac{\partial\beta^{(1)}}{\partial n}\cdot n\frac{\partial\beta^{(2)}}{\partial n}\cdot n\zeta\cdot n d\sigma\\
	-\frac{\mu k}{t^2}\int_{\partial\tilde{\phi}(\Omega)}\frac{\partial w^{(1)}}{\partial n}\frac{\partial w^{(2)}}{\partial n}\zeta\cdot n d\sigma
	+\tilde{\gamma}\int_{\tilde{\phi}(\Omega)}\left(w^{(1)}w^{(2)}+\frac{t^2}{12}\beta^{(1)}\cdot\beta^{(2)}\right)\mathrm{div}\zeta dy.
	\end{multline}	
	
	This, combined with (\ref{maybefirst}), (\ref{second}), (\ref{degei}), concludes the proof.\hfill $\Box$\\

In the case of domain perturbations depending real analytically on one scalar parameter, it is possible to apply the Rellich-Nagy Theorem which allows to conclude that the eigenvalues
splitting from a multiple eigenvalue of multiplicity $m$ are described by $m$ real-analytic functions. Namely, we have the following theorem which can be proved by applying \cite[Cor.~2.28]{lala} combined with Lemma~\ref{tomare}.

\begin{thm}\label{nagy}Let $\Omega$ be a bounded open set in ${\mathbb{R}}^N$ and $t>0$. Let $\tilde \phi \in {\rm BLip }(\Omega)$ and $\{\phi_{\epsilon}\}_{\epsilon\in {\mathbb{R}}} \subset {\rm BLip}(\Omega)$ be a family depending real-analytically on $\epsilon$ such that $\phi_0=\tilde \phi$. Let $\tilde\gamma$ be an eigenvalue of $ R_{\tilde\phi (\Omega),t}$
of multiplcity $m$, with $\tilde \gamma =\gamma_{n,t}[\tilde \phi ]=\dots =\gamma_{n+m-1,t}[\tilde \phi]$ for some $n\in {\mathbb{N}}$. Then there exists an open interval ${\mathcal{I}}$ containing zero and $m$ real-analytic functions $g_{1},\dots , g_m$ from ${\mathcal{I}}$ to ${\mathbb{R}}$ such that $\{\gamma_{n,t}[ \phi_{\epsilon} ],\dots ,\gamma_{n+m-1,t}[\phi_{\epsilon}]\}=\{g_1(\epsilon ), \dots  , g_m(\epsilon) \}$ for all $\epsilon \in {\mathcal{I}}$. Moreover, if $\tilde \phi (\Omega)$  is an open set of class $C^{1,1}$ then the derivatives 
 $g'_1(0), \dots , g_m'(0)$ at zero of the functions $g_1, \dots , g_m$ coincide with the eigenvalues of the matrix $(D_{ij})_{i,j\in \{1, \dots , m\}}$ defined by
	\begin{multline}D_{ij}	=-\frac{\mu}{12}\int_{\partial\tilde{\phi}(\Omega)}\frac{\partial\beta^{(i)}}{\partial n}\cdot\frac{\partial\beta^{(j)}}{\partial n}\zeta\cdot n d\sigma 
-\frac{\mu+\lambda}{12}\int_{\partial\tilde{\phi}(\Omega)}\frac{\partial\beta^{(i)}}{\partial n}\cdot n\frac{\partial\beta^{(j)}}{\partial n}\cdot n\zeta\cdot n d\sigma \\
-\frac{\mu k}{t^2}\int_{\partial\tilde{\phi}(\Omega)}\frac{\partial w^{(i)}}{\partial n}\frac{\partial w^{(j)}}{\partial n}\zeta\cdot n d\sigma ,
	\end{multline}
where $(\beta ^{(i)}, w^{(i)})$, $i=1,\dots ,m$, is an orthonormal basis in ${\mathcal{L}}^2_t(\tilde\phi (\Omega ))$ of the eigenspace associated with $\tilde\gamma$. 
\end{thm}

\section{Isovolumetric perturbations}

Given a bounded open set $\Omega$ in ${\mathbb{R}}^N$,
we consider isovolumetric  domain perturbations, which means that we  consider  transformations $\phi \in {\rm BLip}(\Omega)$ satisfying the volume constraint
\begin{equation}
\label{volume}
| \phi (\Omega )|={\rm constant}.
\end{equation}
  It is then natural to consider the real-valued functional $V$  defined on ${\rm BLip}(\Omega)$  by
\begin{equation}
V[\phi ]={\rm Vol }\, \phi (\Omega ),
\end{equation}
for all $\phi \in {\rm BLip}(\Omega)$.  We recall the following

\begin{defin} Let $\Omega$ be a bounded open set in ${\mathbb{R}}^N$. Let ${\mathcal {F}}$ be a real-valued differentiable map defined
on an open subset of  ${\rm BLip}(\Omega)$. We say that $\tilde \phi \in {\rm BLip}(\Omega)$
is a critical point for ${\mathcal{F}}$ with volume constraint
 if
\begin{equation}
\label{ker}
{\rm ker }\, d_{|_{\phi =\tilde\phi}}V\subset  {\rm ker }\, d_{|_{\phi=\tilde\phi}} {\mathcal{F}}.
\end{equation}
\end{defin}
As is well-known this definition is  related to the problem of finding  local extremal points for the problems 
$$
\min_{V[\phi ]={\rm const}}{\mathcal{F}}[\phi ]\ \ \ {\rm or}\ \ \ \max_{   V[\phi ]={\rm const}} {\mathcal{F}}[\phi].
$$ 
Indeed if $\phi$ is a local minimizer or maximizer of a function ${\mathcal{F}}$ under condition (\ref{volume})   then inclusion (\ref{ker}) holds.

The following theorem can be proved  using formula (\ref{deriv}), by observing that
$d_{|_{\phi=\tilde{\phi}}}V[\psi]=\int_{\partial\tilde{\phi}(\Omega)}(\psi\circ\tilde{\phi}^{(-1)})\cdot n d\sigma$ and by using the Lagrange Multipliers Theorem.  

	\begin{thm}
	\label{moltiplicatori}
	 Let $\Omega$ be a bounded open set in $\mathbb{R}^N$ and $t>0$. Let $F$ be a non-empty finite
	subset of $\mathbb{N}$ and $s\in \{1,\dots , |F|\}$. Let $\tilde{\phi}\in \Theta_{\Omega}[F]$ be such that
	$\tilde{\phi}(\Omega)$ is of class $C^{1,1}$. Then  $\tilde{\phi}$ is a critical point for $\Gamma_{F,t}^{(s)}$  with volume constraint  if and only if there exists an orthonormal basis $(\beta^{(1)}, w^{(1)}),\dots,(\beta^{(|F|)}, w^{(|F|)})$ in ${\mathcal{L}}_t^2(\tilde\phi (\Omega))$ of the eigenspace associated with the eigenvalue
	$\gamma_{F,t}[\tilde{\phi}]$ and there exists $c\in {\mathbb{R}}$ such that
	\begin{equation}
	\label{lacondizione}
	\sum_{l=1}^{|F|}\left(
		\frac{\mu}{12}\left|\frac{\partial\beta^{(l)}}{\partial n}\right|^2
	+\frac{\mu+\lambda}{12}\left(\frac{\partial\beta^{(l)}}{\partial n}\cdot n\right)^2
	+\frac{\mu k}{t^2}\left(\frac{\partial w^{(l)}}{\partial n}\right)^2\right)
	=c \mathrm{\ on\ }\partial\tilde{\phi}(\Omega).
	\end{equation}
	\end{thm}

	As in the case of the Laplace operator discussed  in \cite{lalcri} and  polyharmonic operators considered in \cite{buoso, buosoplates}, it turns out that if $\tilde\phi (\Omega)$ is a ball then condition (\ref{lacondizione}) is satisfied. In order to prove it, we need the following  lemma. Recall that $\beta$ is thought as a row vector.

	\begin{lemm}
	\label{precedente}
	Let $B$ be a ball in $\mathbb{R}^N$ centered at zero,  $t>0$, and let $(\beta , w)$ be an eigenvector of $R_{B ,t}$ in
	$B$ associated with an eigenvalue $\gamma$. Let  $A$  be an orthogonal linear transformation in $\mathbb{R}^N$ and $M$ the corresponding matrix. Then
		also $( (\beta \circ A)M, w\circ A)$ is an eigenvector of $R_{B ,t}$ associated with $\gamma$.
	\end{lemm}
	
	{\bf Proof.}  
	First of all, we note that the rotation invariance of the Laplace operator yields
	$$\Delta(  (\beta \circ A) M)=((\Delta \beta ) \circ A)M,\ \ {\rm and}\ \   \Delta(w\circ A)=(\Delta w)\circ A. $$
	Moreover, by standard calculus we have
	\begin{multline*}
	\mathrm{div}( (\beta\circ A) M)=\mathrm{Tr}\left(  M^T \nabla (\beta\circ A)\right)
	=\mathrm{Tr}\left(M^T((\nabla\beta)\circ A)M\right)=(\mathrm{div}\beta)\circ A ,
	\end{multline*}
where ${\rm Tr}$ denotes the trace of a matrix, 	and
	$$\nabla\mathrm{div}( (\beta\circ A)M)=\nabla((\mathrm{div}\beta)\circ A)=((\nabla\mathrm{div}\beta)\circ A)M.$$
	
By using the previous identities and the fact that  $(\beta ,w )$ is a solution to (\ref{class}), we get 
\begin{eqnarray}
	\label{ruotato}\lefteqn{
	-\frac{\mu}{12}\Delta( (\beta\circ A)M)-\frac{\mu+\lambda}{12}  \nabla\mathrm{div}( (\beta\circ A)M) -\frac{\mu k}{t^2}(\nabla( w\circ A) - (\beta\circ A)M)  }\nonumber \\
& &= -\frac{\mu}{12}((\Delta \beta)\circ A)M-\frac{\mu+\lambda}{12} ((\nabla\mathrm{div}\beta)\circ A)M   -\frac{\mu k}{t^2}(  (\nabla w)\circ A - (\beta\circ A))M \nonumber \\
& & 		=\frac{\gamma t^2}{12}(\beta\circ A )M,
\end{eqnarray}
and 
$$-\frac{\mu k}{t^2}(\Delta (w\circ A)- (\mathrm{div}( (\beta\circ A) M) )=-\frac{\mu k}{t^2}(\Delta w - \mathrm{div}\beta)\circ A=\gamma w\circ A,$$
which show that $( (\beta \circ A)M, w\circ A)$ is an eigenvector of $R_{B ,t}$ associated with $\gamma$.
\hfill $\Box$\\

We now prove the following
		
	\begin{thm}
	\label{lepalle}
		 Let $B$ be the unit ball in $\mathbb{R}^N$ centered at zero, and let $\gamma$ be an eigenvalue of $R_{B,t}$. Let $F$ be the subset of
		$\mathbb{N}$ of indexes $j$ such that $\gamma_{j,t}[B]=\gamma$. Let $(\beta^{(1)}, w^{(1)}),\dots,(\beta^{(|F|)}, w^{(|F|)})$
		be an orthonormal basis in ${\mathcal{L}}_t^2(B)$ of the eigenspace associated with $\gamma$. Then the functions
\begin{equation}\label{lepalle0}
\sum_{l=1}^{|F|}|\beta^{(l)} |^2,\ \ \sum_{l=1}^{|F|}\left|\frac{\partial\beta^{(l)}}{\partial n} \right|^2,\ \ \sum_{l=1}^{|F|}\left|\frac{\partial\beta^{(l)}}{\partial n}\cdot n  \right|^2,\ \ \sum_{l=1}^{|F|}|w^{(l)} |^2,\ \ \sum_{l=1}^{|F|}\left|\frac{\partial w^{(l)}}{\partial n} \right|^2,
\end{equation}
where $n(x)=x/|x|$ for all $x\in \bar B \setminus \{0\}$,  are radial. In particular, there exists $c\in\mathbb{R}$ such that condition (\ref{lacondizione})	holds.
	\end{thm}

	{\bf Proof.}
		Let $O_N(\mathbb{R})$ denote the group of orthogonal linear transformations in $\mathbb{R}^N$, and let $A\in O_N(\mathbb{R})$
		be a transformation with  associated matrix $M$. By Lemma~\ref{precedente}	it follows that $\{ (  (\beta^{(l)}\circ A)M   , w^{(l)}\circ A ):l=1,\dots, |F|\}$ is another orthonormal basis
		of the eigenspace associated with $\gamma$. 
		Since both $\{(\beta^{(l)}, w^{(l)}):l=1,\dots, |F|\}$ and $\{( (\beta^{(l)}\circ A)M , w^{(l)}\circ A):l=1,\dots, |F|\}$
		are orthonormal bases, then there exists $S[A]\in O_{|F|}(\mathbb{R})$ with matrix $(S_{ij}[A])_{i,j=1,\dots,|F|}$ such that

\begin{equation}\label{biorto}		
( (\beta^{(j)}\circ A)M, w^{(j)}\circ A)=\sum_{l=1}^{|F|}S_{jl}[A](\beta^{(l)}, w^{(l)}).
\end{equation}

By (\ref{biorto}) we deduce that 
\begin{equation}
\label{biorto1}
(\beta\circ A)M=S[A]\beta \ \ {\rm and}\ \ w\circ A=S[A]w,
\end{equation}
where $\beta$ denotes the $l\times N$-matrix, the rows of which are given by the row vectors $\beta ^{(j)}$, and $w$ is the column vector the entries of which are given by $w^{(j)}$.

By the first equality in (\ref{biorto1}) we have $(\beta \beta^T)\circ A=S[A]\beta \beta^TS[A]^T$, hence
\begin{equation}\label{orto2}
\sum_{l=1}^{|F|}|\beta^{(l)}\circ A |^2= {\rm Tr}\, [  (\beta \beta^T)\circ A ]={\rm Tr}\, [ S[A]\beta \beta^TS[A]^T]= {\rm Tr}\, [\beta \beta ^T]= 
\sum_{l=1}^{|F|}|\beta^{(l)} |^2.
\end{equation}
 By the arbitrary choice of $A$ we deduce by (\ref{orto2}) that  $\sum_{l=1}^{|F|}|\beta^{(l)} |^2$ is a radial function. Similarly, using the second equality in (\ref{biorto1}), one can prove that $\sum_{l=1}^{|F|}|w^{(l)} |^2$ is a radial function as well. 

We now consider the other functions in (\ref{lepalle0}). By differentiating in the radial direction $n$ the first equality in (\ref{biorto1}), we have that for every $j=1,\dots , l$ and $s=1,\dots , N$,
\begin{equation}
\label{orto3}
\sum_{r,h,k=1}^N\frac{\partial\beta  ^{(j)}_r}{\partial x_h}\circ AM_{hk}M_{rs}n_k=\sum_{l=1}^{|F|}\sum_{k=1}^NS_{jl}[A]\frac{\partial \beta_s^{(l)}}{\partial x_k}n_k.
\end{equation}
Taking into account that $Mn=n\circ A$ we deduce by (\ref{orto3}) that 
\begin{equation}
\label{orto4}
\left(  \frac{\partial \beta }{\partial n}\circ A   \right)M=S[A]\frac{\partial \beta }{\partial n}. 
\end{equation}
By proceeding as in (\ref{orto2}) we get that $\sum_{l=1}^{|F|}\left|\frac{\partial\beta^{(l)}}{\partial n} \right|^2$ is a radial function. 

By multiplying both sides of 
(\ref{orto4}) by $n$ we also get 
\begin{equation}
\left(  \frac{\partial \beta }{\partial n}\cdot n   \right)\circ A=S[A]\frac{\partial \beta }{\partial n}\cdot n,  
\end{equation}
which implies that  $\sum_{l=1}^{|F|}\left|\frac{\partial\beta^{(l)}}{\partial n}\cdot n \right|^2$ is a radial function. Similarly, one can prove that the last function in (\ref{lepalle0})
is radial. \hfill $\Box$\\
		
	 Combining all the results in
	this section we get the following
	
	\begin{thm}
	\label{puntini}  
	Let $\Omega$ be a bounded open set in $\mathbb{R}^N$. Let $\tilde{\phi}\in {\rm BLip}(\Omega)$ be such that
	$\tilde{\phi}(\Omega)$ is a ball. Let $\tilde{\gamma}$ be an eigenvalue of  $R_{\tilde\phi (\Omega), t}$
	and let $F$ be the set of indexes $j\in\mathbb{N}$ such that $\gamma_{j,t}[\tilde{\phi}(\Omega)]=\tilde{\gamma}$.
	Then  for all $s=1,\dots,|F|$ the elementary symmetric function $\Gamma_{F,t}^{(s)}$ has a critical point at $\tilde{\phi}$ with volume constraint.
	\end{thm}

{\bf Acknowledgments.} The authors are very thankful to Professors Carlo Lovadina, Sergei V. Rogosin  and Luis M. Hervella-Nieto for  useful discussions and references.  The authors acknowledge financial support from
the research project `Singular perturbation problems for differential operators',  Progetto di Ateneo of the University of Padova.
The authors are members of the Gruppo Nazionale per l'Analisi Matematica, la Probabilit\`{a} e le loro Applicazioni (GNAMPA) of the
Istituto Nazionale di Alta Matematica (INdAM).

\vspace{1cm}

\noindent Davide Buoso and Pier Domenico Lamberti \\
Dipartimento di Matematica\\
Universit\`{a} degli Studi di Padova\\
Via Trieste 63\\
35121 Padova\\
Italy\\
e-mail:\\ 
dbuoso@math.unipd.it\\
lamberti@math.unipd.it

\end{document}